\documentclass[12pt]{article}
            
\usepackage{amsmath}               
\usepackage{amsfonts}              
\usepackage{amsthm}                
\usepackage[utf8]{inputenc}
\usepackage{leftidx}
\usepackage{tikz-cd}
\usepackage{mathtools}

\newtheorem{tw}{Theorem}[section]

\newtheorem{defi}{Definition}[section]
\newtheorem{prop}{Proposition}[section]
\newtheorem{rem}{Remark}[section]
\newtheorem{cor}{Corollary}[section]

\title{The Fr\"{o}licher-type inequalities of foliations}
\author{Paweł Raźny}

\begin{document}
\maketitle
\section{Introduction}
The purpose of this article is to adapt the Fr\"{o}licher-type inequality, stated and proven for complex and symplectic manifolds in \cite{D2}, to the case of transversely holomorphic and symplectic foliations. These inequalities provide a criterion for checking whether a foliation transversely satisfies the $\partial\bar{\partial}$-lemma and the $dd^{\Lambda}$-lemma (i.e. whether the basic forms of a given foliation satisfy them). These lemmas are linked to such properties as the formality of the basic de Rham complex of a foliation and the transverse hard Lefschetz property. In particular they provide an obstruction to the existence of a transverse K\"{a}hler structure for a given foliation. In the second section we will provide some information concerning the $d'd''$-lemma for a given double complex $(K^{\bullet,\bullet},d',d'')$ and state the main results from \cite{D2}. We will also recall some basic facts and definitions concerning foliations. In the third section we treat the case of transversely holomorphic foliations. We also give a brief review of some properties of the basic Bott-Chern and Aeppli cohomology theories. In Section 4 we prove the symplectic version of the Fr\"{o}licher-type inequality. The final 3 sections of this paper are devoted to the applications of our main theorems. In them we verify the aforementioned lemmas for some simple examples, give the orbifold versions of the Fr\"{o}licher-type inequalities and show that transversely K\"{a}hler foliations satisfy both the $\partial\bar{\partial}$-lemma and the $dd^{\Lambda}$-lemma (or in other words that our main theorems provide an obstruction to the existence of a transversely K\"{a}hler structure).

\section{Preliminaries}
\subsection{$d'd''$-lemma}
Let $(K^{\bullet,\bullet},d',d'')$ be a double cochain complex of modules, and let $D:=d'+d''$ be its total coboundary operator. We say that $K^{\bullet,\bullet}$ satisfies the $d'd''$-lemma iff:
\begin{equation*}
 Ker(d')\cap Ker(d'')\cap Im(D)=Im(d'd'')
\end{equation*}
The first theorem stated here will give us a couple of conditions equivalent to satysfying the $d'd''$-lemma.
\begin{tw}(\cite{Del}, lemma 5.15) Let $K^{\bullet,\bullet}$ be a bounded double complex with notation as above. Then the following conditions are equivalent:
\begin{itemize}
\item $K^{\bullet,\bullet}$ satisfies the $d'd''$-lemma.
\item $Ker(d'')\cap Im(d')=Im(d'd'')$ and $Ker(d')\cap Im(d'')=Im(d'd'')$.
\item $Ker(d')\cap Ker(d'')\cap (Im(d')+Im(d''))=Im(d'd'')$.
\item $Im(d')+Im(d'')+Ker(D)= Ker(d'd'')$.
\item $Im(d'')+Ker(d')= Ker(d'd'')$ and $Im(d')+Ker(d'')=Ker(d'd'')$.
\item $Im(d')+Im(d'')+(Ker(d')\cap Ker(d''))=Ker(d'd'')$.
\end{itemize}
\end{tw}
For any double chain complex there exist two filtrations of the associated total chain complex given by:
\begin{equation*}
'F^p(Tot^k(K^{\bullet,\bullet})):=\bigoplus_{\substack{r+s=k \\ r \geq p}}K^{r,s} \quad ''F^q(Tot^k(K^{\bullet,\bullet})):=\bigoplus_{\substack{r+s=k \\ s \geq q}}K^{r,s}
\end{equation*}
This filtrations induce two spectral sequences $'E^{p,q}_r$ and $''E^{p,q}_r$ in the standard way (cf.\cite{Ei}). This spectral sequences satisfy:
\begin{eqnarray*}
&'E^{p,q}_0=K^{p,q} \quad 'E^{p,q}_1=H^q_{d''}(K^{p,\bullet})\quad 'E^{p,q}_2=H^p_{d'}(H^q_{d''}(K^{\bullet,\bullet}))
\\
&''E^{p,q}_0=K^{p,q} \quad ''E^{p,q}_1=H^q_{d'}(K^{\bullet,p})\quad ''E^{p,q}_2=H^p_{d''}(H^q_{d'}(K^{\bullet,\bullet}))
\end{eqnarray*}
Furthermore, if $K^{\bullet,\bullet}$ is a first (or third) quadrant double complex, then both associated spectral sequences converge to the total cohomology of the complex $K^{\bullet,\bullet}$.
\begin{tw}(\cite{Del}, Proposition 5.17)\label{SS1}
Let $K^{\bullet,\bullet}$ be a bounded double complex of vector space. If the $d'd''$-lemma holds for $K^{\bullet,\bullet}$, then the two spectral sequences associated to $K^{\bullet,\bullet}$ degenerate at the first page.
\end{tw}

\subsection{The Fr\"{o}licher-type inequality for graded vector spaces}
In this subsection we will state the two main theorems from the paper \cite{D2} concerning the Fr\"{o}licher-type inequality for graded and bigraded vector spaces. These theorems will play a crutial role in proving the basic Fr\"{o}licher-type inequalities for transversely holomorphic and transversely symplectic foliations. Let us start with a $\mathbb{Z}^2$-graded $\mathbb{K}$-vector space $V^{\bullet,\bullet}$, endowed with two endomorphisms $\delta_1$ and $\delta_2$ with order $(\delta_{1,1},\delta_{1,2})$ and $(\delta_{2,1},\delta_{2,2})$ respectively. Furthermore, let us assume that this endomorphisms satisfy the condition:
\begin{equation}\label{cc}
\delta_1^2=\delta_2^2=\delta_1\delta_2+\delta_2\delta_1=0
\end{equation}
Let $V^{\bullet}=\bigoplus\limits_{p+q=\bullet}V^{p,q}$ denote the standard $\mathbb{Z}$-graded vector space asociated with $V^{\bullet,\bullet}$. We define the following cohomology of $V^{\bullet,\bullet}$:
\begin{eqnarray*}
& H^{\bullet}_{tot}(V^{\bullet}):=\frac{Ker(\delta_1+\delta_2)}{Im(\delta_1+\delta_2)}
\quad
 H^{\bullet,\bullet}_{\delta_1}(V^{\bullet,\bullet}):=\frac{Ker(\delta_1)}{Im(\delta_1)}
\\
& H^{\bullet,\bullet}_{A}(V^{\bullet,\bullet}):=\frac{Ker(\delta_1\delta_2)}{Im(\delta_1)+Im(\delta_2)}
\quad
 H^{\bullet,\bullet}_{BC}(V^{\bullet,\bullet}):=\frac{Ker(\delta_1)\cap Ker(\delta_2)}{Im(\delta_1\delta_2)}
\end{eqnarray*}
We call them the total, $\delta_1$, Aeppli and Bott-Chern cohomology respectively. One can also define this cohomology groups in the $\mathbb{Z}$-graded case. 
\begin{defi}
Let $(V^{\bullet,\bullet},\delta_1,\delta_2)$  be as above. We say that $V^{\bullet,\bullet}$ satisfies the $\delta_1\delta_2$-lemma iff $Ker(\delta_1)\cap Ker(\delta_2)\cap Im(\delta_1+\delta_2)=Im(\delta_1\delta_2)$.
\end{defi}
It is apparent that this definition coincides with the one in the previous subsection. We can now proceed to state the first of the two main theorems of this subsection:
\begin{tw}\label{M1}
Let $V^{\bullet,\bullet}$ be as above. Suppose that:
\begin{equation*}
dim_{\mathbb{K}}(H_{\delta_i}^{\bullet}(V^{\bullet}))<\infty
\quad  \text{for } i=1,2
\end{equation*}
Then for every $j\in\mathbb{Z}$ the following inequality holds:
\begin{equation*}
dim_{\mathbb{K}}(H^{j}_{BC}(V^{\bullet}))+dim_{\mathbb{K}}(H^{j}_{A}(V^{\bullet}))\geq dim_{\mathbb{K}}(H^{j}_{\delta_1}(V^{\bullet}))+dim_{\mathbb{K}}(H^{j}_{\delta_2}(V^{\bullet}))
\end{equation*}
Furthermore, the equality holds iff $V^{\bullet,\bullet}$ satisfies the $\delta_1\delta_2$-lemma.
\end{tw}
In order to state this theorem in the $\mathbb{Z}$-graded case we need one additional construction. Given a $\mathbb{Z}$-graded vector space $V^{\bullet}$, endowed with two endomorphisms $\delta_1$ and $\delta_2$ of order $|\delta_1|$ and $|\delta_2|$, satisfying (\ref{cc}) we define the associated $\mathbb{Z}^2$-graded vector space by:
\begin{equation*}
Doub^{p,q}(V^{\bullet}):=V^{|\delta_1|p+|\delta_2|q}
\end{equation*}
We can also extend $\delta_1$ and $\delta_2$ to this bigraded vector space in such a way that their orders are $(1,0)$ and $(0,1)$. With this we can now state the second main theorem of this subsection:
\begin{tw}\label{M2}
Let $(V^{\bullet},\delta_1,\delta_2)$ be as above. Furthermore, let $GCD(|\delta_1|,|\delta_2|)=1$ and:
\begin{equation*}
dim_{\mathbb{K}}(H_{\delta_i}^{\bullet}(V^{\bullet}))<\infty
\quad  \text{for } i=1,2
\end{equation*}
Then for every $j\in\mathbb{Z}$ the following inequality holds:
\begin{equation}\label{Fr}
dim_{\mathbb{K}}(H^{j}_{BC}(V^{\bullet}))+dim_{\mathbb{K}}(H^{j}_{A}(V^{\bullet}))\geq dim_{\mathbb{K}}(H^{j}_{\delta_1}(V^{\bullet}))+dim_{\mathbb{K}}(H^{j}_{\delta_2}(V^{\bullet}))
\end{equation}
Moreover, the following two conditions are equivalent:
\begin{enumerate}
\item The equality in (\ref{Fr}) holds and the spectral sequences associated to the double complex $Doub^{\bullet,\bullet}(V^{\bullet})$ degenerate at the first page.
\item $V^{\bullet}$ satisfies the $\delta_1\delta_2$-lemma.
\end{enumerate}
\end{tw}

\subsection{Foliations}
We shall end the preliminary part of this paper with a brief review of some basic facts concerning foliations and transverse structures. The interested reader is referred to \cite{M1} for a more thorough exposition.
\begin{defi} A codimension q foliation $\mathcal{F}$ on a smooth n-manifold M is given by the following data:
\begin{itemize}
\item An open cover $\mathcal{U}:=\{U_i\}_{i\in I}$ of M.
\item A q-dimensional smooth manifold $T_0$.
\item For each $U_i\in\mathcal{U}$ a submersion $f_i: U_i\rightarrow T_0$ with connected fibers (these fibers are called plaques).
\item For all intersections $U_i\cap U_j\neq\emptyset$ a local diffeomorphism $\gamma_{ij}$ of $T_0$ such that $f_j=\gamma_{ij}\circ f_i$
\end{itemize}
The last condition ensures that plaques glue nicely to form a partition of M consisting of submanifolds of M of codimension q. This partition is called a foliation $\mathcal{F}$ of M and the elements of this partition are called leaves of $\mathcal{F}$.
\end{defi}
We call $T=\coprod\limits_{U_i\in\mathcal{U}}f_i(U_i)$ the transverse manifold of $\mathcal{F}$. The local diffeomorphisms $\gamma_{ij}$ generate a pseudogroup $\Gamma$ of transformations on T (called the holonomy pseudogroup).The space of leaves $M\slash\mathcal{F}$ of the foliation $\mathcal{F}$ can be identified with $T\slash\Gamma$.
\begin{defi}
 A smooth form $\omega$ on M is called basic if for any vector field X tangent to the leaves of $\mathcal{F}$ the following equality holds:
\begin{equation*}
i_X\omega=i_Xd\omega=0
\end{equation*}
Basic 0-forms will be called basic functions henceforth.
\end{defi}
Basic forms are in one to one correspondence with $\Gamma$-invariant smooth forms on T. It is clear that $d\omega$ is basic for any basic form $\omega$. Hence, the set of basic forms of $\mathcal{F}$ (denoted $\Omega^{\bullet}(M\slash\mathcal{F})$) is a subcomplex of the de Rham complex of M. We define the basic cohomology of $\mathcal{F}$ to be the cohomology of this subcomplex and denote it by $H^{\bullet}(M\slash\mathcal{F})$. A transverse structure to $\mathcal{F}$ is a $\Gamma$-invariant structure on T. For example:
\begin{defi}
$\mathcal{F}$ is said to be transversely symplectic if T admits a $\Gamma$-invariant closed 2-form $\omega$ of maximal rank. $\omega$ is then called a transverse symplectic form. As we noted earlier $\omega$ corresponds to a closed basic form of rank q on M (also denoted $\omega$).
\end{defi}
\begin{defi}
$\mathcal{F}$ is said to be transversely holomorphic if T admits a complex structure that makes all the $\gamma_{ij}$ holomorphic. This is equivalent to the existence of an almost complex structure $J$ on the normal bundle $N\mathcal{F}:=TM\slash T\mathcal{F}$ (where $T\mathcal{F}$ is the bundle tangent to the leaves) satisfying:
\begin{itemize}
\item $L_XJ=0$ for any vector field X tangent to the leaves.
\item if $Y_1$ and $Y_2$ are sections of the normal bundle then:
\begin{equation*}
 N_J(Y_1,Y_2):=[JY_1,JY_2]-J[Y_1,JY_2]-J[JY_1,Y_2]+J^2[Y_1,Y_2]=0
 \end{equation*}
where $[$ , $]$ is the bracket induced on the sections of the normal bundle.
\end{itemize}
\end{defi}
\begin{rem}
If $\mathcal{F}$ is transversely holomorphic we have the standard decomposition of the space of complex valued forms $\Omega^{\bullet}({M\slash\mathcal{F},\mathbb{C}})$ into forms of type (p,q) and $d$ decomposes into the sum of operators $\partial$ and $\bar{\partial}$ of order (1,0) and (0,1) respectively. Hence, one can define the Dolbeault double complex $(\Omega^{\bullet,\bullet}({M\slash\mathcal{F},\mathbb{C}}),\partial,\bar{\partial})$, the Fr\"{o}licher spectral sequence and the Dolbeault cohomology as in the manifold case (cf.\cite{Wol2}). 
\end{rem}
\begin{defi}
$\mathcal{F}$ is said to be transversely orientable if T is orientable and all the $\gamma_{ij}$ are orientation preserving. This is equivalent to the orientability of $N\mathcal{F}$.
\end{defi}
\begin{defi}
$\mathcal{F}$ is said to be Riemannian if T has a $\Gamma$-invariant Riemannian metric. This is equivalent to the existence of a Riemannian metric g on $N\mathcal{F}$ with $L_Xg=0$ for all vector fields X tangent to the leaves.
\end{defi}
\begin{defi}
$\mathcal{F}$ is said to be transversely parallelizable if there exist q linearly independent $\Gamma$-invariant vector fields.
\end{defi}
\begin{defi}
A foliation is said to be Hermitian if it is both transversely holomorphic and Riemannian.
\end{defi}
Throughout the rest of this section $\mathcal{F}$ will denote a transversely orientable Riemannian foliation on a compact manifold M. Under these assumptions we shall construct a scalar product on the space of basic forms following \cite{E1}. We start with the SO(q)-principal bundle $p:M^{\#}\rightarrow M$ of orthonormal frames transverse to $\mathcal{F}$. The foliation $\mathcal{F}$ lifts to a p dimensional, tansversely parallelizable, Riemannian foliation $\mathcal{F}^{\#}$ on $M^{\#}$. Furthermore, this foliation is SO(q)-invariant (i.e. for any element $a\in SO(q)$ and any leaf L of $\mathcal{F^{\#}}$, $a(L)$ is also a leaf of $\mathcal{F^{\#}}$). There exists a compact manifold W and a fibre bundle $\pi:M^{\#}\rightarrow W$ with fibers equal to the closures of leaves of $\mathcal{F}^{\#}$. The manifold W is called the basic manifold of $\mathcal{F}$. It is apparent that the SO(q)-invariant smooth functions on W and basic functions on M are in one to one correspondence. In particular, for basic k-forms $\alpha$ and $\beta$ the basic function $g_x(\alpha_x,\beta_x)$ induces a SO(q)-invariant function $\Phi (\alpha,\beta)(w)$ on W (where $g_x$ is the scalar product induced on $\bigwedge^k T_x^*M$ by the Riemmanian structure). With this we can define the scalar product on basic forms:
\begin{equation*}
<\alpha,\beta>:=\int_W \Phi(\alpha,\beta)(w)d\mu (w)
\end{equation*}
Where $\mu$ is the meassure associated to the metric on W. The transverse $*$-operator can be defined fiberwise on the orthogonal complements of the spaces tanget to the leaves in the standard way. This construction can be repeated for complex valued basic forms on Hermitian foliations. We shall finish this section by recalling the notion of a transversely elliptic differential operator. We shall restrict our attention to differential operators on basic forms (cf.\cite{E1},\cite{E3}).
\begin{defi}
A basic differential operator of order m is a linear map $D:\Omega^{\bullet}(M\slash\mathcal{F})\rightarrow\Omega^{\bullet}(M\slash\mathcal{F})$ such that in local coordinates $(x_1,...,x_p,y_1,...,y_q)$ (where $x_i$ are leaf-wise coordinates and $y_j$ are transverse ones) it has the form:
\begin{equation*}
D=\sum\limits_{|s|\leq m}a_s(y)\frac{\partial^{|s|}}{\partial^{s_1}y_1...\partial^{s_q}y_q}
\end{equation*}
where $a_s$ are matrices of appropriate size with basic functions as coefficients. A basic differential operator is called transversely elliptic if its principal symbol is an isomorphism at all points of $x\in M$ and all non-zero, transverse, cotangent vectors at x.
\end{defi}
Due to the correspondence between basic forms of $\mathcal{F}$ and $\Gamma$-invariant forms on the transverse manifold T, a basic differential operator induces a $\Gamma$-invariant differential operator on T. Furthermore, transverse ellipticity of a basic differential operator is equivalent to the ellipticity of its $\Gamma$-invariant counterpart (this is apparent since the principal symbol is defined pointwise).
\begin{tw}(cf.\cite{E1},\cite{E3})
Under the above assumptions the kernel of a transversely elliptic differential operator is finitely dimensional.
\end{tw}

\section{Bott-Chern and Aeppli cohomology theories}
\subsection{Basic Bott-Chern cohomology of foliations}
Let M be a manifold of dimension $n=p+2q$, endowed with a Hermitian foliation $\mathcal{F}$ of complex codimension q. Using the basic Dolbeault double complex we can define the basic Bott-Chern cohomology of $\mathcal{F}$:
\begin{eqnarray*}
H^{\bullet,\bullet}_{BC}(M\slash\mathcal{F}):=\frac{Ker(\partial)\cap Ker(\bar{\partial})}{Im(\partial\bar{\partial})}
\end{eqnarray*}
where the operators $\partial$ and $\bar{\partial}$ are defined as the components of order (1,0) and (0,1) of the operator $d$ restricted to the basic forms (as mentioned earlier). Our main goal in this subsection, is to prove the decomposition theorem for basic Bott-Chern cohomology. To that purpose, we define the operator:
\begin{equation*}
\Delta_{BC}:=(\partial\bar{\partial})(\partial\bar{\partial})^*+(\partial\bar{\partial})^*(\partial\bar{\partial})+
(\bar{\partial}^*\partial)(\bar{\partial}^*\partial)^*+(\bar{\partial}^*\partial)^*(\bar{\partial}^*\partial)
+\bar{\partial}^*\bar{\partial}+\partial^*\partial
\end{equation*}
where by $\partial^*$ and $\bar{\partial}^*$, we mean the operators adjoint to $\partial$ and $\bar{\partial}$,  with respect to the Hermitian product, defined by the transverse Hermitian structure.
\begin{prop}
The operator $\Delta_{BC}$ is transversely elliptic and self-adjoint.
\begin{proof}
Being self-adjoint is easy to see using the definition of $\Delta_{BC}$. To prove ellipticity, we will show that the principal symbol of this operator coincides, on the local quotient manifold, with the principal symbol of the manifold version of $\Delta_{BC}$, which is elliptic as proven in \cite{S1}. In \cite{ER} the authors showed that the operator $\delta$ adjoint to $d$ is given by the formula:
\begin{equation*}
\delta\alpha=(-1)^{(n-p)(k+1)+1}*d\alpha-P\kappa\wedge\alpha*
\end{equation*}
where $\alpha$ is a complex valued basic k-form and $P\kappa$ is a basic 1-form dependent on $\mathcal{F}$ (a slightly modified mean curvature). By splitting $P\kappa$ into forms $\kappa_1$ and $\kappa_2$ of type (1,0) and (0,1) respectively, we get the following formula:
\begin{equation*}
\partial^*=(-1)^{(n-p)(k+1)+1}*\partial-\kappa_1\wedge* \quad \bar{\partial}^*=(-1)^{(n-p)(k+1)+1}*\bar{\partial}-\kappa_2\wedge*
\end{equation*}
Since $*\kappa_1\wedge*$ and $*\kappa_2\wedge*$ are 0-order differential operators, they do not contribute to the principal symbol of $\Delta_{BC}$. Hence, the principal symbol of $\Delta_{BC}$ is the same as the principal symbol of:
\begin{equation*}
\partial\bar{\partial}\bar{\partial}'\partial'+
\bar{\partial}'\partial'\partial\bar{\partial}+
\bar{\partial}'\partial\partial'\bar{\partial}+
\partial'\bar{\partial}\bar{\partial}'\partial
\end{equation*}
where:
\begin{equation*}
\partial'=(-1)^{(n-p)(k+1)+1}*\partial* \quad and \quad \bar{\partial}'=(-1)^{(n-p)(k+1)+1}*\bar{\partial}*
\end{equation*}
which in turn has the same principal symbol as the manifold version of $\Delta_{BC}$ on the local quotient manifold.
\end{proof}
\end{prop}

\begin{tw}(Decomposition of the basic Bott-Chern cohomology)
If M is a compact manifold, endowed with a Hermitian foliation $\mathcal{F}$, then we have the following decomposition:
\begin{equation*}
\Omega^{\bullet,\bullet}(M\slash\mathcal{F},\mathbb{C})=Ker(\Delta_{BC})\oplus Im(\partial\bar{\partial})\oplus (Im(\partial^*)+Im(\bar{\partial}^*))
\end{equation*}
In particular,
\begin{equation*}
H^{\bullet,\bullet}_{BC}(M\slash\mathcal{F})\cong Ker(\Delta_{BC})
\end{equation*}
and the dimension of $H^{\bullet,\bullet}_{BC}(M\slash\mathcal{F})$ is finite.
\begin{proof}
Let $<,>$ be the Hermitian product, induced on $\Omega^{\bullet,\bullet}(M\slash\mathcal{F})$ by the Hermitian metric. Then by linearity, and the definition of adjoint operators, we have the following equality, for any $\omega\in\Omega^{\bullet,\bullet}(M\slash\mathcal{F})$
\begin{equation*}
<\omega,\Delta_{BC}\omega>=||\bar{\partial}^*\partial^*\omega||^2+||\partial\bar{\partial}\omega||^2+||\partial^*\bar{\partial}\omega||^2+||\bar{\partial}^*\partial\omega||^2+||\bar{\partial}\omega||^2+||\partial\omega||^2
\end{equation*}
Hence, it is evident that,
\begin{equation*}
\omega\in Ker(\Delta_{BC})\iff (\partial\omega=0,\text{ }\bar{\partial}\omega=0,\text{ }\bar{\partial}^*\partial^*\omega=0)
\end{equation*}
In other words,
\begin{equation*}
Ker(\Delta_{BC})=Ker(\partial)\cap Ker(\bar{\partial})\cap Ker(\bar{\partial}^*\partial^*)
\end{equation*}
By computing it's orthogonal complement, we get the first part of the theorem. To prove the isomorphism, two observation are required. Firstly, $Ker(\Delta_{BC})\subset Ker(\partial)\cap Ker(\bar{\partial})$. Secondly, $Ker(\partial)\cap Ker(\bar{\partial})$ and $(Im(\partial^*)+Im(\bar{\partial}^*))$ have trivial intersection. Finally, $H^{\bullet,\bullet}_{BC}(M\slash\mathcal{F})$ is finite dimensional, because it is isomorphic to the kernel of a self-adjoint transversely elliptic differential operator.
\end{proof}
\end{tw}

\subsection{Basic Aeppli cohomology of foliations}
Let M and $\mathcal{F}$ be as in the previous section. We define the basic Aeppli cohomology of $\mathcal{F}$ to be:
\begin{equation*}
H^{\bullet,\bullet}_{A}(M\slash\mathcal{F}):=\frac{Ker(\partial\bar{\partial})}{Im(\partial)+Im(\bar{\partial})}
\end{equation*}
We define a basic differential operator, needed for the decomposition theorem for the basic Aeppli cohomology of $\mathcal{F}$:
\begin{equation*}
\Delta_A:=\partial\partial^*+\bar{\partial}\bar{\partial}^*+
(\partial\bar{\partial})^*(\partial\bar{\partial})+(\partial\bar{\partial})(\partial\bar{\partial})^*+(\bar{\partial}\partial^*)^*(\bar{\partial}\partial^*)+
(\bar{\partial}\partial^*)(\bar{\partial}\partial^*)^*
\end{equation*}
\begin{prop}
$\Delta_A$ is a self-adjoint, transversely elliptic operator.
\begin{proof}
This proposition is proven in the exact same way as the analogous proposition for $\Delta_{BC}$.
\end{proof}
\end{prop}
\begin{tw}(Decomposition of basic Aeppli cohomology)
Let M be a compact manifold, endowed with a Hermitian foliation $\mathcal{F}$. Then we have the following decomposition:
\begin{equation*}
\Omega^{\bullet,\bullet}(M\slash\mathcal{F},\mathbb{C})=Ker(\Delta_A)\oplus (Im(\partial)+Im(\bar{\partial}))\oplus Im((\partial\bar{\partial})^*)
\end{equation*}
In particular, there is an isomorphism,
\begin{equation*}
H^{\bullet,\bullet}_A(M\slash\mathcal{F})\cong Ker(\Delta_A)
\end{equation*}
and the dimension of $H^{\bullet,\bullet}_A(M\slash\mathcal{F})$ is finite.
\begin{proof}
By calculations similar as in the Bott-Chern case we get the following equality:
\begin{equation*}
Ker(\Delta_A)=Ker(\partial^*)\cap Ker(\bar{\partial}^*)\cap Ker(\partial\bar{\partial})
\end{equation*}
Computing the orthogonal compliment of $Ker(\Delta_A)$, finishes the first part of the proof. The isomorphism is the consequence of two basic facts. Firstly, $Ker(\Delta_A) \subset Ker(\partial\bar{\partial})$. Secondly, $Ker(\partial\bar{\partial})$ and $Im((\partial\bar{\partial})^*)$ have trivial intersection. Finally, $H^{\bullet,\bullet}_{A}(M\slash\mathcal{F})$ is finitely dimensional, because it is isomorphic to the kernel of a self-adjoint, transversely elliptic differential operator.
\end{proof}
\end{tw}
Finally, we will prove a duality theorem for basic Bott-Chern and Aeppli cohomology. However, for the theorem to work, we need an additional condition on our foliation:
\begin{defi} A foliation $\mathcal{F}$ on M is called homologically orientable if $H^{cod\mathcal{F}}(M\slash\mathcal{F})=\mathbb{R}$.
\end{defi}
\begin{rem} The above condition guaranties that the following equalities hold for basic r-forms:
\begin{equation*}
\partial^*=(-1)^{q(r+1)+1}*\partial*
\quad
\bar{\partial}^*=(-1)^{q(r+1)+1}*\bar{\partial}*
\end{equation*}
where $*$ is the transverse $*$-operator. For general foliations this does not have to be true (c.f. \cite{M1}, appendix B, example 2.3 and \cite{E1}).
\end{rem}
\begin{cor}
If M is a compact manifold endowed with a Hermitian, homologically orientable foliation $\mathcal{F}$, then the transverse star operator induces an isomorphism:
\begin{equation*}
H^{p,q}_{BC}(M\slash\mathcal{F})\rightarrow H^{n-p,n-q}_A(M\slash\mathcal{F})
\end{equation*} 
\begin{proof}
From the proofs of decomposition theorems, we know that
\begin{eqnarray*}
u\in Ker(\Delta_{BC})&\iff & \partial u=\bar{\partial} u=(\partial\bar{\partial})^*u=0
\\
&\iff &  \partial^* (*u)=\bar{\partial}^*(*u)= \partial\bar{\partial}(*u)=0
\\
&\iff &(*u)\in Ker(\Delta_A)
\end{eqnarray*}
which proves the duality thanks to the isomorphism described in the decomposition theorems.
\end{proof}
\end{cor}
\begin{rem} The definition of Bott-Chern and Aeppli cohomology theories are valid for transversely holomorphic foliations. However, all the theorems up until now strongly depend on the transverse Hermitian structure.
\end{rem}

\subsection{Basic Fr\"{o}licher-type inequality}
Let us continue with our main result concerning the basic Bott-Chern and Aeppli cohomologies:
\begin{tw}(Basic Fr\"{o}licher-type inequality)
Let M be a manifold of dimension n endowed with a transversely holomorphic foliation $\mathcal{F}$ of complex codimension q. Let us assume that the basic Dolbeault cohomology of $\mathcal{F}$ are finitely dimensional. Then, for every $k\in\mathbb{N}$, the following inequality holds:
\begin{equation*}
\sum\limits_{p+q=k}(dim_{\mathbb{C}}(H^{p,q}_{BC}(M\slash\mathcal{F}))+dim_{\mathbb{C}}(H^{p,q}_A(M\slash\mathcal{F})))\geq 2dim_{\mathbb{C}}(H^k(M\slash\mathcal{F},\mathbb{C}))
\end{equation*}
Furthermore, the equality holds for every $k\in\mathbb{N}$, iff $\mathcal{F}$ satisfies the $\partial\bar{\partial}$-lemma (i.e. its basic Dolbeault double complex satisfies the $\partial\bar{\partial}$-lemma).
\begin{proof}
By applying Theorem \ref{M1} to our case we get the following inequality:
\begin{equation*}
dim_{\mathbb{C}}(H^{j}_{BC}(M\slash\mathcal{F}))+dim_{\mathbb{C}}(H^{j}_{A}(M\slash\mathcal{F}))\geq dim_{\mathbb{C}}(H^{j}_{\partial}(M\slash\mathcal{F}))+dim_{\mathbb{C}}(H^{j}_{\bar{\partial}}(M\slash\mathcal{F}))
\end{equation*}
So all that is left to prove is that the right hand side is bigger than the doubled complex dimension of the basic cohomology of $\mathcal{F}$. Let us consider the spectral sequences associated to the basic Dolbeault double complex. The first page of this spectral sequences are the basic $\partial$ and $\bar{\partial}$ cohomology of $\mathcal{F}$, while their final page in both cases is the basic cohomology. This leads us to the conclusion:
\begin{equation*}
dim_{\mathbb{C}}(\bigoplus\limits_{p+q=k} H^{p,q}_{\bar{\partial}}(M\slash\mathcal{F}))\geq dim_{\mathbb{C}}(H^*(M\slash\mathcal{F},\mathbb{C}))
\end{equation*}
and
\begin{equation*}
dim_{\mathbb{C}}(\bigoplus\limits_{p+q=k} H^{p,q}_{\partial}(M\slash\mathcal{F}))\geq dim_{\mathbb{C}}(H^*(M\slash\mathcal{F},\mathbb{C}))
\end{equation*}
This finishes the proof of the inequality. Now if the equality in the theorem holds then in particular the equality in Theorem \ref{M1} holds, which is equivalent to the ${\partial\bar{\partial}}$-lemma. If on the other hand the ${\partial\bar{\partial}}$-lemma holds then the spectral sequences associated to the basic Dolbeault double complex degenerate at the first page (This is the consequence of Theorem \ref{SS1}). This fact together with Theorem \ref{M1} gives us the desired equality.
\end{proof}
\end{tw}
We will now treat the special case when $\mathcal{F}$ is a Hermitian foliation on a closed manifold M. As it was proven in \cite{E1} the basic Dolbeault cohomology has finite dimension in this case. Hence we get the following corollary:
\begin{cor}
Let $\mathcal{F}$ be a Hermitian foliation on a closed manifold M. Then for all $k\in\mathbb{N}$ the following inequality holds:
\begin{equation*}
\sum\limits_{p+q=k}(dim_{\mathbb{C}}(H^{p,q}_{BC}(M\slash\mathcal{F}))+dim_{\mathbb{C}}(H^{p,q}_A(M\slash\mathcal{F})))\geq 2dim_{\mathbb{C}}(H^k(M\slash\mathcal{F},\mathbb{C}))
\end{equation*}
Furthermore, the equality holds for every $k\in\mathbb{N}$, iff $\mathcal{F}$ satisfies the $\partial\bar{\partial}$-lemma (i.e. it's basic Dolbeault double complex satisfies the $\partial\bar{\partial}$-lemma).
\end{cor}

\section{$dd^{\Lambda}$ and $d+d^\Lambda$ cohomology theories}
\subsection{Basic definitions}
Throughout this section let $\mathcal{F}$ be a transversely symplectic foliation of codimension 2q on a manifold M with basic symplectic form $\omega$. Let us start by defining the symplectic star operator for $\mathcal{F}$. The transverse symplectic form makes the fibers of the normal bundle $N\mathcal{F}$ into a symplectic vector space. Hence, we can define a nondegenerate pairing $\tilde{G}$ on the sections of $\bigwedge^*N^*\mathcal{F}$ in the standard way. We can restrict $\tilde{G}$ to a pairing G on basic forms.
\begin{defi} The symplectic star operator is a linear operator
\begin{equation*}
*_s :\Omega^{k}(M\slash\mathcal{F})\rightarrow\Omega^{2n-k}(M\slash\mathcal{F})
\end{equation*}
uniquely defined by the formula:
\begin{equation*}
\alpha_1\wedge *_s\alpha_2=G(\alpha_1,\alpha_2)\frac{\omega^q}{q!}
\end{equation*}
where $\alpha_1$ and $\alpha_2$ are arbitrary basic k-forms.
\end{defi}
The symplectic $*$-operator is an isomorphism. Using this operator we can define a couple of other important operators on $k$-forms which will be used further in this paper:
\begin{equation*}
L(\alpha):=\omega\wedge\alpha \quad \Lambda(\alpha):=*_s L*_s(\alpha) \quad d^{\Lambda}\alpha:=(-1)^{k+1}*_s d*_s(\alpha)=d\Lambda(\alpha)-\Lambda d(\alpha)
\end{equation*}
The equality in the definition of $d^{\Lambda}$ is a consequence of the exact same local computations as in the manifold case. We can use this operator to define cohomology theories similar to those reviewed in the previous chapter:
\begin{eqnarray*}
H^{\bullet}_{d^{\Lambda}}(M\slash\mathcal{F})&:=&\frac{Ker(d^{\Lambda})}{Im(d^{\Lambda})}
\\ H^{\bullet}_{d+d^{\Lambda}}(M\slash\mathcal{F})&:=&\frac{Ker(d+d^{\Lambda})}{Im(dd^{\Lambda})}
\\ H^{\bullet}_{dd^{\Lambda}}(M\slash\mathcal{F})&:=&\frac{Ker(dd^{\Lambda})}{Im(d)+Im(d^{\Lambda})}
\end{eqnarray*}
It is easy to see that the basic $d^{\Lambda}$-cohomology is simply the basic cohomology with reversed gradation.

\subsection{Symplectic spectral sequence}
We would like to prove the Fr\"{o}licher-type inequality for the newly defined $d+d^{\Lambda}$- and $dd^{\Lambda}$-cohomology theories. However, in order to apply Theorem $\ref{M2}$ we need to prove that the appropriate spectral sequences degenerate at the first page. To do so, we are going to adapt the proof presented in \cite{C1}. This proof is made for complex valued forms, but due to the universal coefficients theorem the degeneration of the spectral sequences at the first page for complex and real valued forms are equivalent. The idea used there, is to change the gradation on the space of forms, show that this change induces a spectral sequence equivalent in some regard to the original one and then prove that this new spectral sequence degenerates at the first page. So first let us state the theorem which will produce the new grading on the space of complex valued basic forms:
\begin{tw}(\cite{C1}, Theorem 2.2) Given a symplectic vector space $(V,\omega)$, the subspaces of  $\wedge^{\bullet}V\otimes\mathbb{C}$ of the form:
\begin{equation*}
U^{q-k}:= \{e^{i\omega}e^{\frac{\Lambda}{2i}}\alpha\text{ }|\text{ } \alpha\in\wedge^{k}V\otimes\mathbb{C}\}
\end{equation*}
form a decomposition. In particular $\wedge^{k}V\otimes\mathbb{C}$ is isomorphic to $U^{q-k}$.
\end{tw}
We can apply this theorem by taking the normal bundle $N\mathcal{F}$. The basic symplectic form will then make each fibre of this bundle into a symplectic vector space. Furthermore, we note that the isomorphism in the theorem applied fiberwise sends basic forms into basic forms since $\omega$ is basic. We shall show that the inverses of $e^{i\omega}$ and $e^{\frac{\Lambda}{2i}}$ also send basic forms to basic forms. Let us assume that $\alpha$ is a non-basic k-form on the normal bundle. The degree k component of $e^{i\omega}\alpha$ is equal to $\alpha$. Which means, that $e^{i\omega}\alpha$ is non-basic as well. This proves that the inverse of $e^{i\omega}$ sends basic forms to basic forms. The proof for $e^{\frac{\Lambda}{2i}}$ is similar. We conclude that $e^{i\omega}e^{\frac{\Lambda}{2i}}$ is an isomorphism on basic forms. Moreover, basic forms decompose as in the previous theorem.
\begin{tw} For any complex valued basic form $\alpha$ the following equality holds:
\begin{equation*}
d(e^{i\omega}e^{\frac{\Lambda}{2i}}\alpha)=e^{i\omega}e^{\frac{\Lambda}{2i}}(d\alpha -\frac{1}{2i}d^{\Lambda}\alpha)
\end{equation*}
\begin{proof}
By the definition of $d^{\Lambda}$  the equality $d^\Lambda=\Lambda d-d^{}\Lambda$ holds. Hence, by induction and the fact that $d^{\Lambda}$ and $\Lambda$ commute, we get the equality:
\begin{equation*}
d\Lambda^{k}=\Lambda^kd+k\Lambda^{k-1}d^{\Lambda}
\end{equation*}
This allows us to make the following computation for any complex valued basic k-form $\alpha$:
\begin{eqnarray*}
d(e^{i\omega}e^{\frac{\Lambda}{2i}}\alpha)=e^{i\omega}d(e^{\frac{\Lambda}{2i}}\alpha)=e^{i\omega}\sum\limits_{j=0}^{\infty}d(\frac{\Lambda^k}{(2i)^kk!}\alpha)=
\\
=e^{i\omega}\sum\limits_{j=0}^{\infty}(\frac{\Lambda^k}{(2i)^kk!}d\alpha-\frac{\Lambda^{k-1}}{(2i)^k(k-1)!}d^{\Lambda}\alpha)=
e^{i\omega}e^{\frac{\Lambda}{2i}}(d\alpha-\frac{1}{2i}d^{\Lambda}\alpha)
\end{eqnarray*}
\end{proof}
\end{tw}
Now let us denote by $\partial$ and $\bar{\partial}$ the projections of $d|_{U^k}$ onto $U^{k+1}$ and $U^{k-1}$ respectively. We state several consequences of the previous theorem which will finish the construction of our alternate spectral sequence:
\begin{eqnarray*}
d&=&\partial+\bar{\partial}
\\
\partial (e^{i\omega}e^{\frac{\Lambda}{2i}}\alpha)&=&-e^{i\omega}e^{\frac{\Lambda}{2i}}\frac{1}{2i}d^{\Lambda}\alpha
\\
\bar{\partial}(e^{i\omega}e^{\frac{\Lambda}{2i}}\alpha)&=&e^{i\omega}e^{\frac{\Lambda}{2i}}d\alpha
\end{eqnarray*}
Having made this preparation we can proceed to the proof of our desired result:
\begin{tw}\label{degen}
Let M be a manifold endowed with a transversely symplectic foliation $\mathcal{F}$ with basic symplectic form $\omega$. Then both spectral sequences of the double complex $Doub^{\bullet,\bullet}(\Omega^{\bullet}(M\slash\mathcal{F}),d,d^{\Lambda})$ degenerate at the first page.
\begin{proof} As we have mentioned before equivalently we can prove that the spectral sequence induced by the double complex associated to complex valued basic forms degenerates at the first page (due to the universal coefficients theorem).  This in turn is equivalent to the degeneration at the first page of the spectral sequences induced by the complex $Doub^{\bullet,\bullet}(U^{\bullet},\partial,\bar{\partial})$ (by the previous theorem). We shall prove that the sequence with $\bar{\partial}$-cohomology on its first page degenerates at the first page (since the proof for the second spectral sequence is identical). Since this sequence is periodic and $d=\partial+\bar{\partial}$ we know that this sequence converges to the basic cohomology of $\mathcal{F}$ (with the even$\backslash$odd gradation). However, using the previous theorem we know that it's first page is isomorphic to the basic cohomology of $\mathcal{F}$. This proves that the spectral sequences of  $Doub^{\bullet,\bullet}(U^{\bullet},\partial,\bar{\partial})$ degenerate at the first page. Hence, the same is true for the sequences induced by $Doub^{\bullet,\bullet}(\Omega^{\bullet}(M\slash\mathcal{F}),d,d^{\Lambda})$.
\end{proof}
\end{tw}

\subsection{Symplectic Fr\"{o}licher-type inequality}
Using the final result from the previous section and Theorem \ref{M2} we obtain the following theorem:
\begin{tw}
Let M be a manifold endowed with a transversely symplectic foliation $\mathcal{F}$ of codimension 2q with transverse symplectic form $\omega$. If the basic cohomology of $\mathcal{F}$ have finite dimension then the following inequality holds for any $j\in\mathbb{N}$:
\begin{equation*}
dim(H^j_{d+d^{\Lambda}}(M\slash\mathcal{F}))+dim(H^j_{dd^{\Lambda}}(M\slash\mathcal{F}))\geq dim(H^{j}(M\slash\mathcal{F}))+dim(H^{2q-j}(M\slash\mathcal{F}))
\end{equation*}
Furthermore, the equality holds if and only if $\mathcal{F}$ satisifies the $dd^{\Lambda}$-lemma (i.e. the complex of basic forms satisfies the $dd^{\Lambda}$-lemma).
\end{tw}
As in the complex case the theorem is greatly simplified in the case of riemannian foliations:
\begin{cor}
Let M be a compact manifold endowed with a transversely symplectic, riemannian, homologically orientable foliation $\mathcal{F}$ of codimension 2q with transverse symplectic form $\omega$. The following inequality holds for any $j\in\mathbb{N}$:
\begin{equation*}
dim(H^j_{d+d^{\Lambda}}(M\slash\mathcal{F}))+dim(H^j_{dd^{\Lambda}}(M\slash\mathcal{F}))\geq 2dim(H^{j}(M\slash\mathcal{F}))
\end{equation*}
Furthermore, the equality holds if and only if $\mathcal{F}$ satisifies the $dd^{\Lambda}$-lemma.
\begin{proof}
In this case we already know that the basic cohomology of $\mathcal{F}$ is finitely dimensional (proven in e.g.\cite{E1}). Since the foliation is homologically orientable, the $j$-th and $(2q-j)$-th basic cohomology are isomorphic (proven in \cite{E2}).
\end{proof}
\end{cor}

\section{Applications to orbifolds}
The results of this section are motivated by the fact that the leaf space of a riemannian foliation with compact leaves on a compact manifold is a compact orbifold (cf.\cite{R1} for basic definitions and properties of orbifolds). Furthermore every orbifold can be realized as a space of leaves of a Riemannian foliation (see \cite{M1} and \cite{Gi}). Moreover, foliated structures on this foliations project onto equivariant structures on the orbifolds of leaves (also basic forms are isomorphic as complexes to differential forms on the orbifold). Due to this and the fact that the proofs from the previous sections can be rewritten here almost word for word we allow ourselves to be a little lax with the notation (unless we prove a result without it's counterpart in the previous sections). It is worth noting that versions of some of the theorems presented in this section were proven in \cite{D4} in a more straightforward manner for global quotient orbifolds.

\subsection{Complex orbifolds}
Let X be a complex orbifold, and let $(\Omega^{\bullet,\bullet}(X),\partial,\bar{\partial})$ be the orbifold Dolbeault double complex, with total coboundary operator $d=\partial+\bar{\partial}$. We define the following cohomology groups:
\begin{eqnarray*}
H^{\bullet,\bullet}_{\partial}(X):=\frac{Ker(\partial)}{Im(\partial)}
\quad
H^{\bullet,\bullet}_{A}(X):=\frac{Ker(\partial\bar{\partial})}{Im(\partial)+Im(\bar{\partial})}
\\
H^{\bullet,\bullet}_{\bar{\partial}}(X):=\frac{Ker(\bar{\partial})}{Im(\bar{\partial})}
\quad
H^{\bullet,\bullet}_{BC}(X):=\frac{Ker(\partial)\cap Ker(\bar{\partial})}{Im(\partial\bar{\partial})}
\end{eqnarray*}
Due to the observation above we get the decomposition theorems and Fr\"{o}licher-type inequality as simple consequences of their foliated counterparts (after remarking that the isomorphism between differential forms on the orbifold and basic forms of the corresponding foliaiton respect $\partial$ and $\bar{\partial}$). We start of again by defining the self-adjoint elliptic operators:
\begin{eqnarray*}
\Delta_{BC}:=(\partial\bar{\partial})(\partial\bar{\partial})^*+(\partial\bar{\partial})^*(\partial\bar{\partial})+
(\bar{\partial}^*\partial)(\bar{\partial}^*\partial)^*+(\bar{\partial}^*\partial)^*(\bar{\partial}^*\partial)
+\bar{\partial}^*\bar{\partial}+\partial^*\partial
\\
\Delta_A:=\partial\partial^*+\bar{\partial}\bar{\partial}^*+
(\partial\bar{\partial})^*(\partial\bar{\partial})+(\partial\bar{\partial})(\partial\bar{\partial})^*+(\bar{\partial}\partial^*)^*(\bar{\partial}\partial^*)+
(\bar{\partial}\partial^*)(\bar{\partial}\partial^*)^*
\end{eqnarray*}
This allows us to state the following two theorems:
\begin{tw}(Decomposition for Bott-Chern cohomology of orbifolds)
If X is a compact complex Hermitian orbifold, then we have the following decomposition:
\begin{equation*}
\Omega^{\bullet,\bullet}(X)=Ker(\Delta_{BC})\oplus Im(\partial\bar{\partial})\oplus (Im(\partial^*)+Im(\bar{\partial}^*))
\end{equation*}
In particular,
\begin{equation*}
H^{\bullet,\bullet}_{BC}(X)\cong Ker(\Delta_{BC})
\end{equation*}
and the dimension of $H^{\bullet,\bullet}_{BC}(X)$ is finite.
\end{tw}
\begin{tw}(Decomposition for Aeppli cohomology of orbifolds)
Let X be a compact complex Hermitian orbifold. Then we have the following decomposition
\begin{equation*}
\Omega^{\bullet,\bullet}(X)=Ker(\Delta_A)\oplus (Im(\partial)+Im(\bar{\partial}))\oplus Im((\partial\bar{\partial})^*)
\end{equation*}
In particular, there is an isomorphism:
\begin{equation*}
H^{\bullet,\bullet}_A(X)\cong Ker(\Delta_A)
\end{equation*}
and the dimension of $H^{\bullet,\bullet}_A(X)$ is finite.
\end{tw}
Finally, we state our main result for complex orbifolds:
\begin{tw}(Fr\"{o}licher-type inequality for orbifolds)
Let X be a compact complex orbifold of complex dimension n. Then, for every $k\in\mathbb{N}$, the following inequality holds:
\begin{equation*}
\sum\limits_{p+q=k}(dim_{\mathbb{C}}(H^{p,q}_{BC}(X))+dim_{\mathbb{C}}(H^{p,q}_A(X)))\geq 2dim_{\mathbb{C}}(H^k_{dR}(X,\mathbb{C}))
\end{equation*}
Furthermore, the equality holds, for every $k\in\mathbb{N}$, iff X satisfies the $\partial\bar{\partial}$-lemma (i.e. it's Dolbeault double complex satisfies the $\partial\bar{\partial}$-lemma).
\end{tw}

\subsection{Symplectic orbifolds}
The focus of this section is on the orbifold version of the symplectic Fr\"{o}licher-type inequality. We start by stating the simple version which is a corollary of it's foliated counterpart. After that, we will relate this theorem to the hard Lefschetz property.
\newline\indent First of all, let us recall that on symplectic orbifolds we also have the symplectic star operator on forms defined in the standard way. Keeping this in mind we can also define the operators L, $\Lambda$, $d^{\Lambda}$ as it is done for manifolds and as we did earlier for foliations. Now we can define cohomology theories that are of interest to us:
\begin{eqnarray*}
H^{\bullet}_{d^{\Lambda}}(X)&:=&\frac{Ker(d^{\Lambda})}{Im(d^{\Lambda})}
\\ H^{\bullet}_{d+d^{\Lambda}}(X)&:=&\frac{Ker(d+d^{\Lambda})}{Im(dd^{\Lambda})}
\\ H^{\bullet}_{dd^{\Lambda}}(X)&:=&\frac{Ker(dd^{\Lambda})}{Im(d)+Im(d^{\Lambda})}
\end{eqnarray*}
\begin{rem} The results from the previous section concerning the degeneration of the spectral sequences can also be adapted to the orbifold case (either by repeating the proof from the foliated case or applying the result of the previous section to a foliation representing the orbifold in question).
\end{rem}
\begin{tw} Let $(X,\omega)$ be a 2n-dimensional symplectic orbifold with finitely dimensional de Rham cohomology. Then for every $k\in\mathbb{N}$ the following inequality holds:
\begin{equation*}
dim(H^k_{d+d^{\Lambda}}(X))+dim(H^k_{dd^{\Lambda}}(X))\geq dim(H^k_{dR}(X))+dim(H^{2n-k}_{dR}(X))
\end{equation*}
Furthermore, the equality holds iff X satisfies the $dd^{\Lambda}$-lemma.
\end{tw}
Our final goal in this section is to prove that an orbifold satisfies the $dd^{\Lambda}$-lemma if and only if it has the hard Lefschetz property:
\begin{tw}
Let X be a symplectic closed orbifold. The following conditions are equivalent:
\begin{enumerate}
\item X has the hard Lefschetz property.
\item The morphism of complexes $(\Omega^{\bullet}(X),d^{\Lambda})\rightarrow (\Omega^{\bullet}(X)\slash d\Omega^{\bullet}(X),d^{\Lambda})$ induces an isomorphism in cohomology.
\item Any class in the de Rham cohomology of X has a representative which is $d^{\Lambda}$-closed.
\item $X$ satisfies the $dd^{\Lambda}$-lemma.
\end{enumerate}
\begin{proof} In \cite{BC} the equivalence between the first and the third condition was proven in the closed case. Injectivity of the induced map in 2 can be equivalently written as $Im(d)\cap Ker(d^{\Lambda})\subset Im(d^{\Lambda})$. After passing to the symplectic complement we get $Ker(d)\subset Im(d)+Ker(d^{\Lambda})$, which in turn is equivalent to condition 3. 3 implies 4 is done word for word as in \cite{merk}. 4 implies $Im(d)\cap Ker(d^{\Lambda})\subset Im(d^{\Lambda})$ by definition. With this the equivalence of 3 and 4 is proven. We finish the proof by noting that 4 implies the surjectivity of the induced map in 2. 
\end{proof}
\end{tw}
We finish this section with two corollaries of the previous theorem.

\begin{cor}
Let M be a closed manifold with a transversely symplectic, Riemannian, homologically orientable foliation $\mathcal{F}$ with closed leaves. The following conditions are equivalent:
\begin{enumerate}
\item $\mathcal{F}$ has the hard Lefschetz property.
\item $\mathcal{F}$ satisfies the $dd^{\Lambda}$-lemma.
\item The morphism of complexes $(\Omega^{\bullet}(M\slash\mathcal{F}),d^{\Lambda})\rightarrow (\Omega^{\bullet}(M\slash\mathcal{F})\slash d\Omega^{\bullet}(X),d^{\Lambda})$ induces an isomorphism in cohomology.
\item Any class in the basic cohomology of $\mathcal{F}$ has a representative which is $d^{\Lambda}$-closed.
\end{enumerate}
\begin{proof} Since the leaf space of such a foliation is a closed symplectic orbifold this is a consequence of the two previous theorems.
\end{proof}
\end{cor}
\begin{rem}
The compact leaves assumption can be dropped by using the same method as we use to prove the $dd^{\Lambda}$-lemma in the transversely K\"{a}hler case. All that needs to be done is to show that under our assumptions there exists a compatible transverse almost complex structure. But the second proof of Proposition 2.50 (contractibility of the space of compatible almost complex structures) in \cite{Du} is easily addaptable to the transverse setting under the given assumptions. Existence of a compatible transverse almost complex structure is a simple corollary. Note however that the existence of some inner product is required for the proof to work and hence our foliation has to be Riemannian in order to apply this line of reasoning.
\end{rem}
\begin{cor}
 Let $(X,\omega)$ be a 2n-dimensional closed symplectic orbifold. Then for every $k\in\mathbb{N}$ the following inequality holds:
\begin{equation*}
dim(H^k_{d+d^{\Lambda}}(X))+dim(H^k_{dd^{\Lambda}}(X))\geq 2dim(H^k_{dR}(X))
\end{equation*}
Furthermore, the equality holds iff X has the hard Lefschetz property.
\end{cor}

\section{Examples}
In this section we are going to verify the $dd^{\Lambda}$-lemma and the $\partial\bar{\partial}$-lemma for some of the foliations presented in \cite{Wol}. Keeping this in mind we will compute just enough cohomology groups to prove or disprove the aforementioned lemmas.
\subsection{A transversely symplectic foliation not satisfying the $dd^{\Lambda}$-lemma}
Let N be the Lie group of real matrices of the form:
\begin{equation*}
\begin{bmatrix}
 1 & x & t & z 
\\ 0 & 1 & 0 & y
\\ 0 & 0 & 1 & 0
\\ 0 & 0 & 0 & 1 
\end{bmatrix}
\end{equation*}
To simplify the notation we shall denote such a matrix by (x,y,z,t) in this subsection. We fix an irrational number s and define a subgroup of N (denoted by $\Gamma$) of matrices of the form:
\begin{equation*}
\begin{bmatrix}
 1 & x_1+ sx_2 & t & z_1 + sz_2
\\ 0 & 1 & 0 & y
\\ 0 & 0 & 1 & 0
\\ 0 & 0 & 0 & 1 
\end{bmatrix}
\end{equation*}
where $x_1, x_2. y, z_1, z_2, t$ are integers. We consider the left action of $\Gamma$ on N. Let us consider a third Lie group $U:=(\mathbb{R}^6,\square)$, with the group operation:
\begin{eqnarray*}
&(a_1,a_2,b,c_1,c_2,d)\square(x_1,x_2,y,z_1,z_2,t)=\\
&=(a_1+x_1,a_2+x_2,b+y,c_1+z_1+a_1y,c_2+z_2+a_2y,d+t)
\end{eqnarray*}
The group $\Gamma$ is isomorphic to $(\mathbb{Z}^6,\square)$. Furthernore, there is a submersion $u:U\rightarrow N$ given by the formula:
\begin{equation*}
u(x_1,x_2,y,z_1,z_2,t)\rightarrow (x_1+sx_2,y,z_1+sz_2,t)
\end{equation*}
The foliation on U defined by this submersion is $(\mathbb{Z}^6,\square)$-equivariant. Hence, it descends to a foliation $\mathcal{F}$ on $U\slash (\mathbb{Z}^6,\square)$. The basic forms of $\mathcal{F}$ correspond to $\Gamma$-invariant forms on N. Furthermore, this foliation is transversely symplectic due to the invariant, closed, nondegenerate form:
\begin{equation*}
\omega:= dx\wedge (dz-xdy)+ dy\wedge dt
\end{equation*} 
Let us start by making some observations:
\begin{rem} The orbits of the group $\Gamma$ are dense in the x and z directions and have period one in the other two directions. Hence, basic functions of this foliation coincide with the smooth functions on a torus (the values depend on y and t).
\end{rem}
\begin{rem}
By taking $dx,dy,dz-xdy,dt$ as the orthonormal basis we define a transverse Riemannian metric on our foliation. This means that all the basic cohomology groups are finitely dimensional and hence we can apply the Fr\"{o}licher-type inequality to determine if this foliation satisfies the $dd^{\Lambda}$-lemma. We rename the chosen basis of one forms to $\alpha_1,\alpha_2,\alpha_3,\alpha_4$.
\end{rem}
Since the equality in the appropriate version of the Fr\"{o}licher-type inequality fails for second basic cohomology, we will focus our attention on 2-forms. Let:
\begin{equation*}
\alpha:=\sum\limits_{i<j\leq 4} f_{ij}\alpha_{i}\alpha_{j}
\end{equation*}
be an arbitrary 2-form. By straightforward computation one can see that in the basic, $d+d^{\Lambda}$ and $dd^{\Lambda}$ cohomology, the parts $f_{13}\alpha_1\alpha_3$ and $f_{24}\alpha_2\alpha_4$ both generate a copy of $\mathbb{R}$ in cohomology and have no other influence on the cohomology. Hence, in further computation we can omit them. Without these parts the vector spaces $Ker(d)$ and $Ker(d)\cap Ker(d^{\Lambda})$ become equal. However, by computing the appropriate images we can see that $\alpha_1\alpha_2$ belongs to the image of $d$ and does not belong to the image of $dd^{\Lambda}$. This means that the dimension of the second $(d+d^{\Lambda})$-cohomology is greater than the dimension of the second basic cohomology. The argument is similar for the $dd^{\Lambda}$-cohomology. Here the images can be proven equal (modulo the part generated by $f_{13}\alpha_1\alpha_3$ and $f_{24}\alpha_2\alpha_4$). However, the form $\alpha_3\alpha_4$ belongs to the kernel of $dd^{\Lambda}$ and doesn't belong to the kernel of d. This proves the inequality:
\begin{equation*}
dim(H^2_{d+d^{\Lambda}}(M\slash\mathcal{F}))+dim(H^2_{dd^{\Lambda}}(M\slash\mathcal{F}))> 2dim(H^{2}(M\slash\mathcal{F}))
\end{equation*}
In particular, the $dd^{\Lambda}$-lemma does not hold for this foliation.

\subsection{A transversely holomorphic foliation not satisfying the $\partial\bar{\partial}$-lemma}
We present now our second example. Let N be the Lie group of upper-triangular matrices in $GL(3,\mathbb{C})$. And let $\Gamma$ be its subgroup consisting of the matrices of the form:
\begin{equation*}
\begin{bmatrix}
 1 & z_1+sz'_1 & z_3+sz'_3
\\ 0 & 1 & z_2
\\ 0 & 0 & 1 
\end{bmatrix}
\end{equation*}
where $z_i,z'_i$ are Gauss integers (denoted $\mathbb{Z}[i]$) and s is a fixed irrational number. We again consider the left action of $\Gamma$ on N and take the quotient with respect to this action. As in the previous example we consider another Lie group $U:=(\mathbb{C}^5,\square)$ with group operation:
\begin{equation*}
(u_1,...,u_5)\square (z_1,...,z_2)=(z_1+u_1,z_2+u_2,z_3+u_3,z_4+u_4+u_1z_3,z_5+u_5+u_2z_3)
\end{equation*}
It is evident that $\Gamma$ is isomorphic to $((\mathbb{Z}[i])^5,\square)$ and as before there is a submersion $u:U\rightarrow N$ given by:
\begin{equation*}
u(z_1,...,z_5)=(z_1+sz_2,z_3,z_4+sz_5)
\end{equation*}
This submersion defines a foliation $\mathcal{F}$ on $M:=U\slash ((\mathbb{Z}[i])^5,\square)$.
\begin{rem}
By choosing the invariant orthonormal basis of 1-forms:
\begin{equation*}
dz_1, dz_2, dz_3-z_1dz_2, d\bar{z_1}, d\bar{z_2}, d\bar{z_3}-\bar{z_1}d\bar{z_2}
\end{equation*}
we define a transverse riemannian metric on this foliation. It is easy to compute that this foliation is homologically orientable. Hence, we can use the Fr\"{o}licher-type inequality to determine whether this foliation satisfies the $\partial\bar{\partial}$-lemma
\end{rem}
\begin{rem} As in the previous example the basic function of this foliation depend only on the $z_2$ variable and are the same as the smooth functions on a torus.
\end{rem}
The equality in the Fr\"{o}licher-type inequality fails for the first basic cohomology. By making some straightforward computation we get:
\begin{eqnarray*}
&dim_{\mathbb{C}}(H^1(M\slash\mathcal{F},\mathbb{C}))=4 \\
&dim_{\mathbb{C}}(H^{1,0}_{BC}(M\slash\mathcal{F}))=dim_{\mathbb{C}}(H^{0,1}_{BC}(M\slash\mathcal{F}))=2 \\
&dim_{\mathbb{C}}(H^{1,0}_A(M\slash\mathcal{F}))=dim_{\mathbb{C}}(H^{0,1}_A(M\slash\mathcal{F}))=3
\end{eqnarray*}
which proves that this foliation does not satisfy the $\partial\bar{\partial}$-lemma.

\subsection{A transversely holomorphic foliation satysfying the $\partial\bar{\partial}$-lemma}
In our final example let N be the Lie group consisting of complex matrices of the form:
\begin{equation*}
\begin{bmatrix}
 1 & \bar{z_1} & z_2
\\ 0 & 1 & z_1
\\ 0 & 0 & 1 
\end{bmatrix}
\end{equation*}
and let $\Gamma$ be it's subgroup consisting of matrices of the form:
\begin{equation*}
\begin{bmatrix}
 1 & \bar{z_1}+s\bar{z_1}' & z_2+sz'_2+s^2z''_2
\\ 0 & 1 & z_1+sz'_1
\\ 0 & 0 & 1 
\end{bmatrix}
\end{equation*}
where $z_1,z'_1,z_2,z_2',z_2''$ are Gauss integers and s is a fixed irrational number. Let's consider the group $U:=(\mathbb{C}^5,\square)$ with group operation:
\begin{eqnarray*}
&(u_1,...,u_5)\square (z_1,...z_5)=\\
&=(u_1+z_1,u_2+z_2,u_3+z_3+\bar{u}_1z_1,u_4+z_4+\bar{u}_1z_2+\bar{u}_2z_1,u_5+z_5+\bar{u}_2z_2)
\end{eqnarray*}
It is clear that $\Gamma$ is isomorphic to $((\mathbb{Z}[i])^5,\square)$ and as before there is a $((\mathbb{Z}[i])^5,\square)$-equivariant submersion $u:U\rightarrow N$ given by:
\begin{equation*}
u(z_1,...,z_5)=(z_1+sz_2,z_3+sz_4+s^2z_5)
\end{equation*}
This submersion defines a non-K\'{a}hler foliation $\mathcal{F}$ on $M:=U\slash ((\mathbb{Z}[i])^5,\square)$ (cf.\cite{Wol}). It is easy to see that the only basic functions on this foliated manifold are constant. By taking the invariant basis of forms $(dz_1,dz_2+\bar{z_1}dz_{1}, d\bar{z_1}, d\bar{z_2}+z_1d\bar{z_1})$ one can compute that all the complex Fr\"{o}licher-type equalities hold and hence this foliation satisfies the $\partial\bar{\partial}$-lemma.
\begin{rem}
This foliation also has a transversely symplectic structure but the symplectic Fr\"{o}licher-type equality fails in the second cohomology.
\end{rem}
\begin{rem}
Analogously as in the manifold case the $\partial\bar{\partial}$-lemma implies the formality of the complex of basic forms for transversely holomorphic foliations. In particular, this foliation is formal. This implication is proven exactly as in the manifold case (cf.\cite{Del}).
\end{rem}

\section{Transversely K\"{a}hler foliations}
In this section we will prove the $\partial\bar{\partial}$-lemma and $dd^{\Lambda}$-lemma for transversely K\"{a}hler foliations. As a consequence of this the Fr\"{o}licher-type equalities provide an obstruction, which is relatively easy to check, to the existence of a transversely K\"{a}hler structure on a given foliation. We will start with the $\partial\bar{\partial}$-lemma. The following theorem was formulated in \cite{Wol2}.
\begin{tw}
Let $\mathcal{F}$ be a homologically orientable transversely K\"{a}hler foliation on a compact manifold M. Then $\mathcal{F}$ satisfies the $\partial\bar{\partial}$-lemma. In particular $\mathcal{F}$ is formal.
\begin{proof} Let $\alpha\in\Omega^{p,q}(M\slash\mathcal{F})$ be a $\partial$-closed, $\bar{\partial}$-closed and d-exact form. We define the 3 transversely elliptic operators:
\begin{eqnarray*}
\Delta:=dd^*+d^*d
\\
\Delta_{\partial}:=\partial\partial^*+\partial^*\partial
\\
\Delta_{\bar{\partial}}:=\bar{\partial}\bar{\partial}^*+\bar{\partial}^*\bar{\partial}
\end{eqnarray*}
With the help of this operators one can state and prove the Hodge decomposition for basic and basic Dolbeault cohomology. Furthermore, in the transversely K\"{a}hler case the kernels of this operators are pairwise equal (cf.\cite{E1}). Since $\alpha$ is $d$-exact, it is orthogonal to the space of $\Delta$-harmonic forms. Hence, it is also orthogonal to the space of $\Delta_{\partial}$-harmonic forms. Since $\alpha$ is also $\partial$-closed, it has to be $\partial$-exact as well (by the Hodge decomposition for basic Dolbeault cohomology). Let $\beta$ be a basic form such that $\alpha=\partial\beta$. By applying the Hodge decomposition again we get $\beta=h+\bar{\partial}\eta+\bar{\partial}^*\xi$, where $h$ is $\Delta_{\bar{\partial}}$-harmonic ($h$ is also $\Delta_{\partial}$-harmonic). It suffices to prove that $\bar{\partial}^*\xi$ is $\Delta_{\partial}$-harmonic. Since $\alpha$ is $\bar{\partial}$-closed, we get $\bar{\partial}\bar{\partial}^*\partial\xi=0$ (we use here also one of the K\"{a}hler identities proven in \cite{E1}). Using the scalar product associated to the transverse riemannian metric we get:
\begin{equation*}
||\bar{\partial}^*\partial\xi||^2=<\partial\xi,\bar{\partial}\bar{\partial}^*\partial\xi>=0
\end{equation*}
Which means that $\alpha=\partial\bar{\partial}\eta$.
\end{proof}
\end{tw}
Let us now prove the $dd^{\Lambda}$-lemma for transversely K\"{a}hler homologically orientable foliations:
\begin{tw} Let $\mathcal{F}$ be a transversely K\"{a}hler homologically orientable foliation on a compact manifold M. Then the following statements are true:
\begin{enumerate}
\item $\mathcal{F}$ has the hard Lefschetz property.
\item $\mathcal{F}$ satisfies the $dd^{\Lambda}$-lemma.
\item The morphism of complexes $(\Omega^{\bullet}(M\slash\mathcal{F}),d^{\Lambda})\rightarrow (\Omega^{\bullet}(M\slash\mathcal{F})\slash d\Omega^{\bullet}(X),d^{\Lambda})$ induces an isomorphism in cohomology.
\item Any class in the basic cohomology of $\mathcal{F}$ has a representative which is $d^{\Lambda}$-closed.
\end{enumerate}
\begin{proof} The first condition is satisfied under our assumptions (cf.\cite{E1}). Due to our foliation being transversely K\"{a}hler, $d^{\Lambda}$ is in fact the symplectic adjoint of $d$. Indeed, in this case the symplectic form can be written in terms of the transverse metric and the complex structure in the standard way. Furthermore the symplectic star operator can be written in terms of the Hodge star operator and the complex structure in the standard way. Hence, for any k-form $\beta$ and $(k-1)$-form $\alpha$ we have:
\begin{eqnarray*}
\omega(d\alpha,\beta)&=&<Jd\alpha,\beta>=<\alpha,(-1)^{k+1}*d*J\beta>
\\&=&<J\alpha,(-1)^{k+1}J*dJ*\beta>=\omega(\alpha,d^{\Lambda}\beta)
\end{eqnarray*}
where J is an automorphism induced on the basic forms, by the transverse complex structure and $\omega$ is a non-degenerate pairing induced by the transverse symplectic structure. Keeping this in mind one can repeat the proves from the symplectic orbifold section to deduce the 3 latter conditions from the hard Lefschetz property.
\end{proof}
\end{tw}


\begin{thebibliography}{99}
\bibitem{R1} A. Adem,J. Leida ,Y. Ruan.
\emph{Orbifolds and Stringy Topology}.
Cambridge Press, (2007)
\bibitem{D4} D. Angella,
\emph{Cohomologies of certain orbifolds}. J. Geom.Phys., 71, 117–126 (2013)
\bibitem{D1} D. Angella,
\emph{Cohomological aspects in complex non-K\"{a}hler geometry}. Springer (2014)
\bibitem{D2} D. Angella, A. Tomassini, 
\emph{Inequalities à la Frölicher and cohomological decompositions}.  J. Noncommut. Geom. 9 no. 2, 505–542 (2015)
\bibitem{D3} D. Angella, A. Tomassini,
\emph{On the $\partial\bar{\partial}$-Lemma and Bott-Chern cohomology}. Invent.
Math. 192(1), 71–81 (2013)
\bibitem{BC} L. Bak, A. Czarnecki,
\emph{A remark on the Brylinski conjecture for orbifolds}. J. AUST MATH SOC vol. 91 (01), 1-12 (2011)
\bibitem{Ei} H. Cartan, S. Eilenberg
\emph{Homological Algebra} Princeton University Press (1956)
\bibitem{C1} G.R. Cavalcanti,
\emph{The decomposition of forms and cohomology of generalized complex manifolds}. J. Geom. Phys. 57(1), 121–132 (2006)
\bibitem{Wol} A. Cordero, R. Wolak
\emph{Examples of foliations with foliated geometric structures}. Pacific J. Math vol. 142 (02), 265-276 (1990)
\bibitem{Wol2} A. Cordero, R. Wolak
\emph{Properties of the basic cohomology of transversely K\"{a}hler foliations}. Rendiconti del Circolo Matematico di Palermo, Serie II 40, 177-188 (1991) 
\bibitem{Del} P. Deligne, Ph.A. Griffiths, J. Morgan, D.P. Sullivan,
\emph{Real homotopy theory of K\"{a}hler manifolds}. Invent. Math. 29(3), 245–274 (1975)
\bibitem{ER} P. Efton, K. Richardson,
\emph{The basic Laplacian of a Riemannian foliation} American J. Math. Vol 118 (6), 1249-1275 (1996)
\bibitem{E1} A. El Kacimi-Alaoui,
\emph{Op\'{e}rateurs transversalement elliptiques sur un feuilletage riemannien et applications}. Compositio Mathematica, 73, 57-106 (1990)
\bibitem{E3} A. El Kacimi-Alaoui,
\emph{Towards a basic index theory}. Proceedings of the Summer School and Workshop
Dirac Operator : Yesterday and Today CAMS-AUB, Beirut 2001, 251-261 (2005)
\bibitem{E2} A. El Kacimi-Alaoui, G. Hector, 
\emph{D\'{e}composition de Hodge basique pour un feuilletage Riemannien}. Ann. Inst. Fourier de Grenoble 36 (3), 207-227 (1986)
\bibitem{Gi}  J. Girbau, A. Haefliger, D. Sundararaman,
\emph{On deformations of transversely holomorphic foliations}. J. Reine Angew. Math., 345,122–147 (1983)
\bibitem{Ma} O. Mathieu,
\emph{Harmonic cohomology classes of symplectic manifolds,
Comment}. Comment. Math. Helvetici 70 (1995) 1-9.
\bibitem{Du} D. McDuff, D. Salamon
\emph{Introduction to symplectic topology}. Clarendon Press Oxford (1998)
\bibitem{merk} S.A. Merkulov
\emph{Formality of canonical symplectic complexes and Frobenius
manifolds}.Int. Math. Res. Not. 14, 727–733 (1998)
\bibitem{M1} P. Molino,
\emph{Riemannian foliations}.  Birkh\"{a}user, 1986. Translated by G. Cairns
\bibitem{S1} M. Schweitzer,
\emph{Autour de la cohomologie de Bott-Chern}.  arXiv:0709.3528 [math.AG], 2007 
\end{thebibliography}
\end{document}